\documentclass{amsart}
\usepackage{bbm}
\usepackage{comment}
\usepackage{amsmath, amsthm, amsfonts, amssymb, graphicx, xcolor}
\usepackage{mathrsfs}
\usepackage{enumitem}
\usepackage{tikz}
\usetikzlibrary{positioning}
\usetikzlibrary{arrows.meta}
\usepackage{caption}
\usepackage[colorlinks,
linkcolor=blue,
anchorcolor=blue,
citecolor=blue]{hyperref}

\makeatletter
\@namedef{subjclassname@2020}{\textup{2020} Mathematics Subject Classification}
\makeatother

\newtheorem{thm}{Theorem}[section]
\newtheorem{lem}{Lemma}[section]
\newtheorem{prop}{Proposition}[section]

\newtheorem{dfn}{Definition}[section]

\begin{document}

\title[Weighted Inequalities for the Fractional Maximal Operator on Trees]{Weighted Inequalities for the Fractional Maximal Operator on Trees with Flow Measures}
\author[S. Lu and F. Zhao]{Shengzhe Lu and Fayou Zhao$^*$}
\address[Shengzhe Lu]{Department of Mathematics, Shanghai University, Shanghai 200444, P.R. China}
\email{lushengzhe@shu.edu.cn}

\address[Fayou Zhao]{Newtouch Center for Mathematics of Shanghai University, Department of Mathematics,
	Shanghai University, Shanghai 200444, P.R. China}
\email{fyzhao@shu.edu.cn}
\thanks{$\ast$ Corresponding author}
\subjclass[2020]{05C05, 42B25, 42B35.}

\begin{abstract}
This paper develops a weighted theory for the fractional Hardy-Littlewood maximal operator $M_{\mu,\gamma}$, $0<\gamma<1$, on infinite trees equipped with flow measures. Such measures are typically nondoubling and exhibit exponential volume growth, so the classical Muckenhoupt theory does not apply directly. By replacing Euclidean balls with admissible trapezoids, a class of fractional $A(p,q)(\mu)$ weights is introduced, and a Muckenhoupt-type theorem is established: for $1\le p<1/\gamma$ and $1/p-1/q=\gamma$, the operator $M_{\mu,\gamma}$ is bounded from $\ell^p_\mu(\omega^p)$ to $\ell^q_\mu(\omega^q)$ if and only if $\omega\in A(p,q)(\mu)$. The endpoint case $p=1$ is characterized via a weak-type estimate. 
\end{abstract}

\keywords{fractional maximal operator; weight function; admissible trapezoids; flow measure}
\maketitle

\section{Introduction}
For a locally integrable function $f$ on $\mathbb{R}^n$ and for $0 \le \alpha < n$, the fractional Hardy-Littlewood maximal operator $\mathcal{M}_\alpha$ is defined by
\[
\mathcal{M}_\alpha f(x) = \sup_{r > 0} \frac{1}{|B(x,r)|^{1-\alpha/n}} \int_{B(x,r)} |f(y)| \, dy, \qquad x \in \mathbb{R}^n,
\]
where $|B(x,r)|$ denotes the $n$-dimensional Lebesgue measure of the ball $B(x,r)$. In the particular case $\alpha = 0$, the fractional maximal operator $\mathcal{M}_\alpha$ reduces to the classical Hardy-Littlewood maximal operator $\mathcal{M}$.

The classical Hardy-Littlewood maximal operator $\mathcal{M}$ satisfies the weak type $(1,1)$ estimate and the strong type $(p,p)$ estimates for every $1<p\le\infty$; this result is known as the Hardy-Littlewood-Wiener theorem \cite{HL1930, WN1939}. A natural and fundamental problem is to characterize the weights $\omega$ on $\mathbb{R}^n$ for which the classical Hardy-Littlewood maximal operator $\mathcal{M}$ is bounded on the weighted Lebesgue space $L^p(\omega):=L^p(\mathbb{R}^n,\omega\,dx)$, $1<p<\infty$. Equivalently, one asks to determine all weights $\omega$ for which there exists a constant $C>0$, depending only on $n$, $p$, and $\omega$, such that
\[
\int_{\mathbb{R}^n} (\mathcal{M}f(x))^p \omega(x)\,dx
\le C \int_{\mathbb{R}^n} |f(x)|^p \omega(x)\,dx
\]
for every $f\in L^p(\omega)\cap L^1_{\mathrm{loc}}(\mathbb{R}^n)$. The endpoint case $p=1$ is governed by the weak type $(1,1)$ inequality
\[
\omega\bigl(\{x\in\mathbb{R}^n:\mathcal{M}f(x)>\lambda\}\bigr)
\le \frac{C}{\lambda}\int_{\mathbb{R}^n}|f(x)|\omega(x)\,dx,
\qquad \lambda>0,
\]
for every $f\in L^1(\omega)\cap L^1_{\mathrm{loc}}(\mathbb{R}^n)$.
This problem is not only of intrinsic interest but also has important applications, since the Hardy-Littlewood maximal function naturally dominates error terms in convergence estimates. Consequently, such weighted inequalities are indispensable for proving weighted mean convergence theorems for orthogonal series expansions; see \cite{M1,M2,MS}.

This problem was completely solved by Muckenhoupt \cite{MB1972}, who introduced the now-classical $A_p$ weight classes. More precisely, a weight $\omega$ (a nonnegative locally integrable function) belongs to the class $A_p$ for $1<p<\infty$ if
\[
\sup_{Q}\bigg(\frac{1}{|Q|}\int_{Q}\omega(x)\,dx\bigg)\bigg(\frac{1}{|Q|}\int_{Q}\omega(x)^{-\frac{1}{p-1}}\,dx\bigg)^{\!p-1}<\infty,
\]
and $\omega\in A_1$ if
\[
\sup_{Q}\bigg(\frac{1}{|Q|}\int_{Q}\omega(x)\,dx\bigg)\Bigl(\operatorname{ess\,sup}_{Q}\omega(x)^{-1}\Bigr)<\infty,
\]
where the supremum is taken over all cubes $Q\subset\mathbb{R}^n$. Muckenhoupt proved that, for $1<p<\infty$, $\omega\in A_p$ is both necessary and sufficient for the boundedness of $\mathcal{M}$ on $L^p(\omega)$. At the endpoint $p=1$, the class $A_1$ is both necessary and sufficient for the weak type $(1,1)$ boundedness of $\mathcal{M}$ with respect to $\omega$.

For the fractional counterpart, the corresponding weighted theory was initiated by Muckenhoupt and Wheeden \cite{MW1974}. Given $0<\alpha<n$, $1<p<n/\alpha$, and $1/q=1/p-\alpha/n$, they introduced the $A(p,q)$ classes and proved that $\omega\in A(p,q)$ is necessary and sufficient for the boundedness of the fractional Riesz potential operator $I_\alpha$ from $L^p(\omega^p)$ to $L^q(\omega^q)$; moreover, the same condition characterizes the corresponding boundedness of the fractional Hardy-Littlewood maximal operator $\mathcal{M}_\alpha$. More precisely, for $1<p<\infty$ and $q\ge 1$, a weight $\omega$ belongs to $A(p,q)$ if
\[
\sup_{Q}\bigg(\frac{1}{|Q|}\int_{Q}\omega(x)^{q}\,dx\bigg)^{\!1/q}
\bigg(\frac{1}{|Q|}\int_{Q}\omega(x)^{-p'}\,dx\bigg)^{\!1/p'}<\infty,
\]
and $\omega\in A(1,q)$ if
\[
\sup_{Q}\bigg(\frac{1}{|Q|}\int_{Q}\omega(x)^{q}\,dx\bigg)^{\!1/q}
\Bigl(\operatorname{ess\,sup}_{Q}\omega(x)^{-1}\Bigr)<\infty,
\]
where $p'$ denotes the conjugate exponent of $p$ (i.e., $1/p+1/p'=1$). Here the supremum is taken over all cubes $Q\subset\mathbb{R}^n$. For $p>1$, the corresponding weighted inequality is of strong type, while the endpoint case $p=1$ is understood in the weak type sense.

The classical weighted theory described above relies crucially on the Euclidean structure and on the doubling property of the Lebesgue measure. It is therefore natural to ask to what extent this theory persists in settings where such structural assumptions are absent. Trees provide a natural and tractable framework for investigating this question. Let $T=(V,E)$ be a tree, where $V$ is the vertex set and $E$ is the edge set, equipped with the counting measure $v$ and the geodesic distance $d$. In this setting, Euclidean balls are replaced by the family $\mathcal{R}$ of admissible trapezoids, and the Lebesgue measure is replaced by a flow measure $\mu$, that is, a positive function on $T$ satisfying the flow condition
\[
\mu(x)=\sum_{y\in s(x)}\mu(y),
\]
where $s(x)$ denotes the set of successors of $x$. The resulting space $(T,d,\mu)$ need not be doubling. Recall that $(T,d,v)$ is doubling if there exists a constant $C>0$ such that
\begin{align*}
	\frac{v(B(x,2r))}{v(B(x,r))}\leq C
	\qquad \text{for all } x\in V,\ r>0,
\end{align*}
where $B(x,r):=\{y\in V:d(x,y)<r\}$ denotes the ball centered at $x$ with radius $r$, and $C$ is independent of $x$ and $r$.

Since the doubling condition plays a central role in Muckenhoupt's characterization of $A_p$ weights on spaces of homogeneous type, the classical Euclidean weighted theory cannot be applied directly to trees. This motivates the development of a weighted theory for maximal operators on trees, based on suitable substitutes for Euclidean balls and the usual $A_p$ conditions.  

To present the relevant results, some definitions and notation are needed. Let $\partial T$ denote the boundary of $T$, namely the set of all points at infinity obtained as endpoints of half-infinite geodesics in $T$. Fix a reference point $\xi\in\partial T$ and regard $T$ as hanging down from $\xi$. This induces a partial order $\leq$ on $T$ defined by
\[
x\leq y \quad \text{if and only if} \quad y\in [x,\xi),
\]
where $[x,\xi)$ denotes the semi-infinite geodesic starting at $x$ and tending to $\xi$. Fix a reference vertex $o\in V$, and let $\{x_j\}_{j\ge 0}$ be an enumeration of $[o,\xi)$ such that $x_0=o$ and $d(x_j,x_k)=|j-k|$ for all $j,k\ge 0$. For every $x\in V$, the level function is defined by
\[
l(x)=\lim_{j\to\infty}(j-d(x,x_j)).
\]
For every $x\in V$, let $s(x)$ denote the set of successors of $x$,
\[
s(x):=\{y\sim x:\ l(y)=l(x)-1\},
\]
and let $p(x)$ denote the set of predecessors of $x$,
\[
p(x):=\{y\sim x:\ l(y)=l(x)+1\}.
\]
A positive function $\mu$ on $V$ is called a flow measure if it satisfies the flow condition
\[
\mu(x)=\sum_{y\in s(x)}\mu(y).
\]
Consider an infinite tree $T$ with root $\xi$ at infinity, and assume that $\mu$ is a flow measure on $T$ (possibly not locally doubling) and that $\omega$ is a weight on $T$.

The introduction of admissible trapezoids is motivated by the following difficulty. Flow measures on trees are typically nondoubling and may grow at least exponentially; consequently, the Euclidean balls fail to satisfy the doubling condition with respect to such measures. Since the classical Calder\'{o}n--Zygmund theory, as well as the theory of Hardy and BMO spaces, relies crucially on the doubling property of balls, a different family of sets is needed. In \cite{MFAM2023}, Levi,  Santagati,  Tabacco and Vallarino introduced the family $\mathcal{R}$ of admissible trapezoids as substitutes for Euclidean balls. 
\begin{dfn}
A trapezoid with root $x\in V$ and heights $h_1,h_2\in\mathbb N$, $h_1<h_2$, is defined by
\[
R_{h_1}^{h_2}(x):=\{y\in V:\ y\leq x,\ h_1\leq d(x,y)<h_2\}.
\]
Fix $\beta\geq 12$. A trapezoid $R=R_{h_1}^{h_2}(x)$ is called admissible if $2\leq h_2/h_1\leq \beta$. Its envelope is defined as
\[
\widetilde{R}:=R_{\lceil h_1/\beta\rceil}^{h_2\beta}.
\]

Let $\mathcal R$ denote the family of all admissible trapezoids and singletons. For every $R\in\mathcal R$, denote by $x_R$ the root of $R$, so that $R=R_{h_1}^{h_2}(x_R)$ for some heights $h_1<h_2$. See Figure~\ref{G1} for an illustration.
\end{dfn}
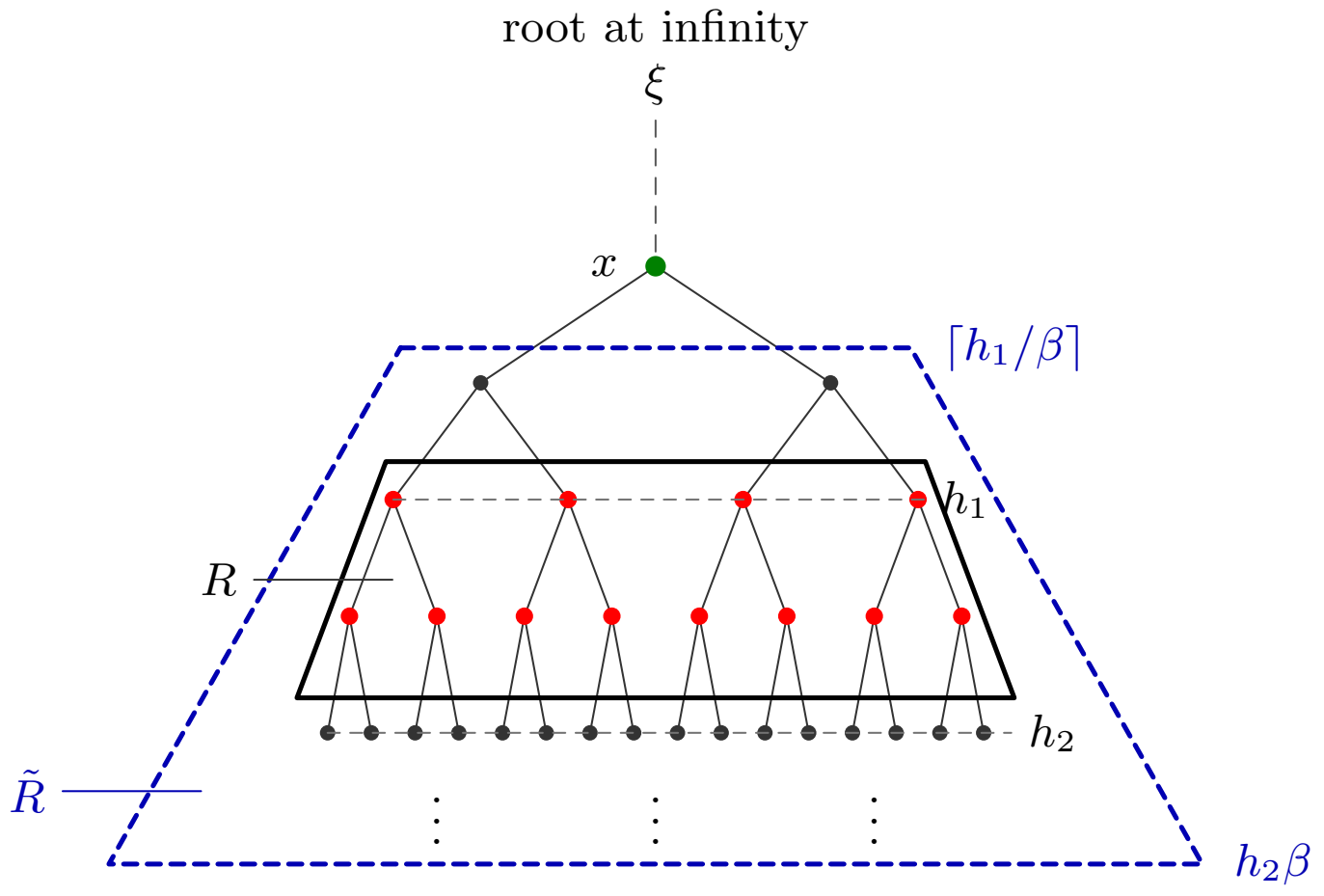
\begin{figure}[htbp]
	\centering
\def\sc{1.00}    
\pgfmathsetlengthmacro\lwTree {0.38pt*\sc}   
\pgfmathsetlengthmacro\lwTrap {0.95pt*\sc}   
\pgfmathsetlengthmacro\aOn    {3.2pt*\sc}    
\pgfmathsetlengthmacro\aOff   {2.4pt*\sc}
\pgfmathsetlengthmacro\eOn    {2.6pt*\sc}    
\pgfmathsetlengthmacro\eOff   {2.0pt*\sc}
\pgfmathsetlengthmacro\lOn    {2.4pt*\sc}    
\pgfmathsetlengthmacro\lOff   {2.0pt*\sc}
\pgfmathsetmacro\fsz{9*\sc}
\pgfmathsetmacro\fsl{11*\sc}
\begin{tikzpicture}[scale=\sc, transform shape,
	line cap=round, line join=round,
	tree/.style   = {draw=black!80, line width=\lwTree},
	lev/.style    = {draw=black!55, line width=\lwTree, dash pattern=on \lOn off \lOff},
	vert/.style   = {circle, fill=black!80, draw=none, inner sep=0pt, minimum size=3.0pt},
	sing/.style   = {circle, fill=red,            draw=none, inner sep=0pt, minimum size=3.4pt},
	sgl/.style    = {circle, fill=green!50!black, draw=none, inner sep=0pt, minimum size=4.0pt},
	axis/.style   = {draw=black!65, line width=\lwTree, dash pattern=on \aOn off \aOff},
	env/.style    = {draw=blue!70!black, line width=\lwTrap, dash pattern=on \eOn off \eOff},
	trap/.style   = {draw=black, line width=\lwTrap, fill=none},      
	lab/.style    = {inner sep=1.5pt, font=\fontsize{\fsz}{\fsl}\selectfont},
	envlab/.style = {lab, text=blue!70!black},
	sglab/.style  = {lab, text=green!50!black},   
	]
	\def\hs{0.30}    
	\def\dy{0.80}    
	\def\Nlv{4}      
	\def\khA{2}      
	\def\khB{4}      
	\def\yhA{-1.34}\def\whA{1.85}   
	\def\yhB{-2.96}\def\whB{2.46}   
	\def\yEnvT{-0.56}\def\whT{1.75} 
	\def\yEnvB{-4.10}\def\whE{3.75} 
	\def\yXi{1.00}   
	
	\draw[axis] (0,\yXi) -- (0,0);
	\node[lab, anchor=south] at (0,\yXi+0.05) {$\xi$};
	\node[lab] at (0,\yXi+0.62) {root at infinity};
	
	\foreach \k/\nb in {0/1,1/2,2/4,3/8}{
		\pgfmathtruncatemacro{\nbm}{\nb-1}
		\foreach \j in {0,...,\nbm}{
			\draw[tree]
			({(2*\j+1)*2^(\Nlv-\k-1)*\hs-2^(\Nlv-1)*\hs},{-\k*\dy})
			--({(4*\j+1)*2^(\Nlv-\k-2)*\hs-2^(\Nlv-1)*\hs},{-(\k+1)*\dy});
			\draw[tree]
			({(2*\j+1)*2^(\Nlv-\k-1)*\hs-2^(\Nlv-1)*\hs},{-\k*\dy})
			--({(4*\j+3)*2^(\Nlv-\k-2)*\hs-2^(\Nlv-1)*\hs},{-(\k+1)*\dy});
		}
	}
	\foreach \k/\nb in {1/2,2/4,3/8,4/16}{
		\pgfmathtruncatemacro{\nbm}{\nb-1}
		\foreach \j in {0,...,\nbm}{
			\ifnum\k>\numexpr\khA-1\relax
			\ifnum\k<\khB
			\node[sing] at ({(2*\j+1)*2^(\Nlv-\k-1)*\hs-2^(\Nlv-1)*\hs},{-\k*\dy}) {};
			\else
			\node[vert] at ({(2*\j+1)*2^(\Nlv-\k-1)*\hs-2^(\Nlv-1)*\hs},{-\k*\dy}) {};
			\fi
			\else
			\node[vert] at ({(2*\j+1)*2^(\Nlv-\k-1)*\hs-2^(\Nlv-1)*\hs},{-\k*\dy}) {};
			\fi
		}
	}
	\foreach \xx in {-1.50,0,1.50}{ \node[lab] at (\xx,-3.70) {$\vdots$}; }
	
	\draw[env] (-\whT,\yEnvT) -- (\whT,\yEnvT) -- (\whE,\yEnvB) -- (-\whE,\yEnvB) -- cycle;
	\node[envlab, anchor=west] at (\whT+0.16,\yEnvT) {$\lceil h_1/\beta\rceil$};
	\node[envlab, anchor=west] at (\whE+0.16,\yEnvB) {$h_2\beta$};
	
	\draw[trap] (-\whA,\yhA) -- (\whA,\yhA) -- (\whB,\yhB) -- (-\whB,\yhB) -- cycle;
	
	\draw[lev] (-1.80,-2*\dy) -- (1.84,-2*\dy);
	\node[lab, anchor=west] at (1.90,-2*\dy) {$h_1$};
	\draw[lev] (-2.25,-4*\dy) -- (2.44,-4*\dy);
	\node[lab, anchor=west] at (2.50,-4*\dy) {$h_2$};
	
	\pgfmathsetmacro\ymid{(\yhA+\yhB)/2}
	\pgfmathsetmacro\whSideA{\whA+(\whB-\whA)*(\yhA-\ymid)/(\yhA-\yhB)}
	\draw[tree] (-\whSideA-0.60,\ymid) -- (-\whSideA+0.35,\ymid);
	\node[lab, anchor=east] at (-\whSideA-0.66,\ymid) {$R$};
	
	\def\yRlab{-3.60}
	\pgfmathsetmacro\whRside{\whT+(\whE-\whT)*(\yEnvT-\yRlab)/(\yEnvT-\yEnvB)}
	\draw[draw=blue!70!black, line width=\lwTree] (-\whRside-0.60,\yRlab) -- (-\whRside+0.35,\yRlab);
	\node[envlab, anchor=east] at (-\whRside-0.66,\yRlab) {$\tilde R$};
	
	\node[sgl] at (0,0) {};
	\node[lab, anchor=east] at (-0.20,0) {$x$};
	
\end{tikzpicture}
\caption{An admissible trapezoid $R_{h_1}^{h_2}(x)$ (black), its envelope $\widetilde{R}$ (blue dotted).}
\label{G1}
\end{figure}

These sets enjoy several essential properties. First, their flow measure is proportional to their height. Second, they admit dyadic decomposition and expansion algorithms. Third, they satisfy the doubling-type estimate
\[
\frac{1}{\widetilde{C}}\mu(R')\le \mu(R)\le \widetilde{C}\mu(Q)
\]
for every admissible trapezoid $R$, its expansion $R'$, and every admissible trapezoid $Q$ in the first decomposition of $R$, where $\widetilde{C}>1$ depends only on $\beta$ and on the bounded-degree constant $c$ of $T$, that is,
\[
|s(x)|\le c \qquad \text{for every } x\in V.
\]
Moreover, by the definition of the envelope and the admissibility condition, there exists a constant $C_\beta>0$, depending only on $\beta$, such that
\[
\mu(\widetilde{R})\le C_\beta \mu(R)
\qquad\text{for every }R\in\mathcal R.
\]
Since every admissible trapezoid has finite height, it is a finite set, and hence
\[
\mu(R)<\infty \qquad \text{for every } R\in\mathcal R.
\]

In view of these properties, admissible trapezoids provide a suitable framework for developing a Calder\'{o}n--Zygmund theory, Hardy spaces, and BMO spaces on trees with nondoubling flow measures; see \cite{MFAM2023} for further details. Building on this framework, the following definition of the uncentred fractional Hardy-Littlewood maximal operator on $(T,d,\mu)$ is introduced in this paper.
\begin{dfn}
Let $0\le \gamma<1$, and let $\mu$ be a flow measure on $T$. For a function $f:V\to\mathbb{C}$ satisfying
\[
\sum_{y\in R}|f(y)|\mu(y)<\infty
\qquad \text{for every } R\in\mathcal R,
\]
the uncentred fractional Hardy-Littlewood maximal operator $M_{\mu,\gamma}$ is defined by
\[
M_{\mu,\gamma}f(x)
:=\sup_{R\in\mathcal R,\,x\in R}
\frac{1}{\mu(R)^{1-\gamma}}
\sum_{y\in R}|f(y)|\mu(y),
\qquad x\in V,
\]
where
\[
\mu(R):=\sum_{y\in R}\mu(y).
\]
The value $M_{\mu,\gamma}f(x)$ is understood in the extended real sense. In the case $\gamma=0$, $M_{\mu,\gamma}$ reduces to the classical Hardy-Littlewood maximal operator $M_\mu$.
\end{dfn}

On trees, the study of the Hardy-Littlewood maximal operator $M_\mu$ began with the unweighted theory. For homogeneous trees, Cowling, Meda and Setti \cite{CMS2010} gave a direct geometric proof of the weak type $(1,1)$ estimate, while Naor and Tao \cite{NT2010} independently established the weak type $(1,1)$ and strong type $(p,p)$ bounds on infinite rooted $k$-ary trees by a combinatorial expansion argument. In the setting of trees with flow measures, Levi, Santagati, Tabacco and Vallarino \cite{MFAM2023} introduced the family $\mathcal R$ of admissible trapezoids on $(T,d,\mu)$ and proved the weak type $(1,1)$ and strong type $(p,p)$ boundedness of $M_\mu$; see \cite[Theorem 3.3]{MFAM2023}. Subsequently, Ottazzi, Santagati and Vallarino \cite{AFM2024} developed a weighted theory in this framework: they introduced the $A_p(\mu)$ weight class and established a Muckenhoupt-type characterization for the boundedness of $M_\mu$ on $L^p(\mu w)$.

For the infinite rooted $k$-ary tree, Ombrosi and Rivera-R\'{i}os \cite{OR2022} studied weighted $L^p$ estimates for $M$ and obtained sufficient conditions for weak and strong type $(p,p)$ bounds. They also proved a Fefferman--Stein inequality and showed that the classical ball $A_p$ condition is not necessary in this setting, revealing a pathological feature absent in the Euclidean case. Very recently, Yang, Yuan and Zhang \cite{YYZ2026} introduced a two-layer Muckenhoupt weight class $\mathcal A_p$ and characterized the boundedness of $M$ on $L^p(w)$ for all $p\in(0,\infty)$ on infinite rooted $k$-ary trees. Their characterization yields, among other things, the optimal quantitative estimate with exponent $1/p$ of the weight constant and a global Sawyer-type testing condition.

Motivated by these developments, we introduce the $A(p,q)(\mu)$ class on $(T,d,\mu)$ for $p\in[1,\infty)$, with admissible trapezoids replacing Euclidean balls. We then investigate the weighted boundedness of the fractional Hardy--Littlewood maximal operator on $(T,d,\mu)$. Specifically, we prove that the fractional Hardy-Littlewood maximal operator $M_{\mu,\gamma}$, associated with the family of admissible trapezoids and the flow measure $\mu$, is bounded from $\ell^p_{\mu}(\omega^p)$ to $\ell^q_{\mu}(\omega^q)$ whenever $0<\gamma<1$, $1\le p<1/\gamma$, and $1/p-1/q=\gamma$. The case $\gamma=0$ reduces to the classical Hardy-Littlewood maximal operator and recovers the corresponding result in \cite{AFM2024}.

\section{The weight class $A(p,q)(\mu)$ and its properties}

This section begins with the definition of the $A(p,q)(\mu)$ weight class on trees, which is motivated by the Euclidean $A(p,q)$ classes recalled in the introduction. In the Euclidean setting, the supremum is taken over all cubes and the Lebesgue measure is used. In the present framework, cubes are replaced by admissible trapezoids and the Lebesgue measure by the flow measure $\mu$. Since $\mu$ is not assumed to be doubling, this definition cannot be obtained by a direct transplantation of the Euclidean one, and it is tailored to the family $\mathcal{R}$ of admissible trapezoids. The resulting class will be used to characterize the weighted boundedness of the fractional Hardy-Littlewood maximal operator $M_{\mu,\gamma}$.

For later use, the following local weighted sequence spaces are introduced. For $R\in\mathcal R$, $1\le p<\infty$, and a weight $\omega$, define
\[
\ell^p_{R,\mu}(\omega)
:=
\bigg\{
f:R\to\mathbb C:
\|f\|_{\ell^p_{R,\mu}(\omega)}
:=
\bigg(\sum_{x\in R}|f(x)|^p\omega(x)\mu(x)\bigg)^{1/p}
<\infty
\bigg\}.
\]
For $R=V$, the notation $\ell^p_{\mu}(\omega)$ is used in place of $\ell^p_{V,\mu}(\omega)$; that is
\[
\ell^p_{\mu}(\omega)
:=
\ell^p_{V,\mu}(\omega)
=
\bigg\{
f:V\to\mathbb C:
\|f\|_{\ell^p_{\mu}(\omega)}
:=
\bigg(\sum_{x\in V}|f(x)|^p\omega(x)\mu(x)\bigg)^{1/p}
<\infty
\bigg\}.
\]

\begin{dfn}\label{dfn2-1}
Let $1<p<\infty$ and $1\le q<\infty$. A weight $\omega$ is called an $A(p,q)(\mu)$ weight if
\[
[\omega]_{A(p,q)(\mu)}
:=\sup_{R\in\mathcal{R}}
\bigg(\frac{1}{\mu(R)}\sum_{y\in R}\omega(y)^q\mu(y)\bigg)^{1/q}
\bigg(\frac{1}{\mu(R)}\sum_{y\in R}\omega(y)^{-p'}\mu(y)\bigg)^{1/p'}
<\infty.
\]
A weight $\omega$ is called an $A(1,q)(\mu)$ weight if
\[
[\omega]_{A(1,q)(\mu)}
:=\sup_{R\in\mathcal{R}}
\bigg(\frac{1}{\mu(R)}\sum_{y\in R}\omega(y)^q\mu(y)\bigg)^{1/q}
\bigg(\operatorname*{ess\,sup}_{y\in R}\omega(y)^{-1}\bigg)
<\infty,
\]
where and in what follows, $p'$ denotes the conjugate exponent of $p$, that is, $1/p+1/p'=1$.
\end{dfn}

For comparison, the Muckenhoupt $A_p(\mu)$ weight classes on $(T,d,\mu)$, introduced by Ottazzi, Santagati, and Vallarino in \cite{AFM2024}, are recalled. For $1<p<\infty$, a weight $\omega$ is called an $A_p(\mu)$ weight if
\[
[\omega]_{A_p(\mu)}
:=
\sup_{R\in\mathcal{R}}
\bigg(\frac{1}{\mu(R)}\sum_{y\in R}\omega(y)\mu(y)\bigg)
\bigg(\frac{1}{\mu(R)}\sum_{y\in R}\omega(y)^{-\frac{1}{p-1}}\mu(y)\bigg)^{p-1}
<\infty.
\]
A weight $\omega$ is called an $A_1(\mu)$ weight if
\[
[\omega]_{A_1(\mu)}
:=
\sup_{R\in\mathcal{R}}
\bigg(\frac{1}{\mu(R)}\sum_{y\in R}\omega(y)\mu(y)\bigg)
\Bigl(\operatorname*{ess\,sup}_{y\in R}\omega(y)^{-1}\Bigr)
<\infty.
\]

 In order to prove the main theorem, some basic properties of the $A(p,q)(\mu)$ weights are needed. These properties clarify the relationship between $A(p,q)(\mu)$ and the classical $A_p(\mu)$ weights and provide equivalent characterizations that will be used in the proof. They are collected in the following proposition.

\begin{prop}\label{prop 2-1}
Let $\omega\in A(p,q)(\mu)$ with $0\leq \gamma<1$, $1\leq p<1/\gamma$, and ${1}/{p}-{1}/{q}=\gamma$.

{\rm(i)} If $p>1$, then $\omega\in A(p,q)(\mu)$ is equivalent to $\omega^q\in A_{1+\frac{q}{p'}}(\mu)$, which in turn is equivalent to $\omega^{-p'}\in A_{1+\frac{p'}{q}}(\mu)$.

{\rm(ii)} If $p>1$ and $\omega\in A(p,q)(\mu)$, then
$\omega^q\in A_q(\mu)$ and $\omega^p\in A_p(\mu)$.

{\rm(iii)} If $p=1$, then
$\omega\in A(1,q)(\mu)$ if and only if $\omega^q\in A_1(\mu)$.

{\rm(iv)} If $\omega\in A(1,q)(\mu)$, then
\[
\lim_{p\to1^{+}}[\omega]_{A(p,q)(\mu)}=[\omega]_{A(1,q)(\mu)}.
\]
	
{\rm(v)} An alternative characterization of $A(1,q)(\mu)$ weights is given as follows:
\[
[\omega]_{A(1,q)(\mu)}
=
\sup_{R\in\mathcal R}
\sup_{\substack{f\in\ell^1_{R,\mu}(\omega^q)\\ \|f\|_{\ell^1_{R,\mu}(\omega^q)}>0}}
\Bigg(
\frac{\frac{1}{\mu(R)}\sum_{x\in R}|f(x)|\mu(x)}
{\frac{1}{\omega_{q,\mu}(R)}\sum_{x\in R}|f(x)|\omega(x)^q\mu(x)}
\Bigg)^{1/q},
\]
where and in what follows,
\[
\omega_{q,\mu}(R):=\sum_{x\in R}\omega(x)^q\mu(x)
\] denotes the total mass of $\omega^q$ on $R$ with respect to $\mu$.

Moreover, for $1<p<1/\gamma$,
\[
[\omega]_{A(p,q)(\mu)}
\gtrsim
[\omega^q]_{A_q(\mu)}^{1/q}
=
\sup_{R\in\mathcal R}
\sup_{\substack{f\in\ell^q_{R,\mu}(\omega^q)\\ \|f\|_{\ell^q_{R,\mu}(\omega^q)}>0}}
\Bigg(
\frac{\left(\frac{1}{\mu(R)}\sum_{x\in R}|f(x)|\mu(x)\right)^q}
{\frac{1}{\omega_{q,\mu}(R)}\sum_{x\in R}|f(x)|^q\omega(x)^q\mu(x)}
\Bigg)^{1/q}.
\]
	
{\rm(vi)} For every $\omega\in A(p,q)(\mu)$ with $1<p<1/\gamma$, one has
\[
[\omega]_{A(p,q)(\mu)}\ge 1,
\]
with equality if and only if $\omega$ is constant $\mu$-a.e. The same conclusion holds for $\omega\in A(1,q)(\mu)$.

{\rm(vii)} If $p=1$, then for all $R_1,R_2\in\mathcal R$ satisfying
\[
\mu(R_1)\le C_0\mu(R_1\cap R_2)
\]
for some constant $C_0>0$, one has
\[
\omega_{q,\mu}(R_1)\le C_0 [\omega]^q_{A(1,q)(\mu)}\omega_{q,\mu}(R_2).
\]
If $1<p<1/\gamma$, the same conclusion holds with $A(p,q)(\mu)$ in place of $A(1,q)(\mu)$:
\[
\omega_{q,\mu}(R_1)\le C_{p,q,C_0}[\omega]^q_{A(p,q)(\mu)}\omega_{q,\mu}(R_2).
\]
\end{prop}

\begin{proof}
The proofs of (i)--(iii) follow an approach similar to that in \cite[Proposition 7.1.5]{G2014}.

{\rm(iv)} Let $q_1=\frac{1}{1-\gamma}$. For $p>1$, it follows from
${1}/{p}-{1}/{q}=\gamma $ that  $q\to q_1$ as $p\to1^+$. By property (i), for $p>1$,
$
[\omega]_{A(p,q)(\mu)}
=
[\omega^q]_{A_{1+q/p'}(\mu)}^{1/q}.
$
Set
$
r(p)=1+{q}/{p'}.
$
Then $r(p)\to1^+$ as $p\to1^+$. Since $\omega\in A(1,q_1)(\mu)$, property (iii) gives $\omega^{q_1}\in A_1(\mu)$ and
\[
[\omega]_{A(1,q_1)(\mu)}
=
[\omega^{q_1}]_{A_1(\mu)}^{1/q_1}.
\]

We first prove the upper bound. For classical $A_r$ weights, if $v\in A_1$, then
\[
\lim_{r\to1^+}[v]_{A_r}=[v]_{A_1},
\]
and the convergence is uniform under the $A_1$ constant; see \cite[Proposition 7.1.5 (7)]{G2014}. Taking $v=\omega^{q_1}$ and using $q\to q_1$, $r(p)\to1^+$, we obtain
\[
\limsup_{p\to1^+}[\omega]_{A(p,q)(\mu)}
=
\limsup_{p\to1^+}[\omega^q]_{A_{r(p)}(\mu)}^{1/q}
\le
[\omega^{q_1}]_{A_1(\mu)}^{1/q_1}
=
[\omega]_{A(1,q_1)(\mu)}.
\]

For the lower bound, fix $R\in\mathcal R$ and set
\[
F_R(p)
:=
\bigg(\frac{1}{\mu(R)}\sum_{x\in R}\omega(x)^q\mu(x)\bigg)^{1/q}
\bigg(\frac{1}{\mu(R)}\sum_{x\in R}\omega(x)^{-p'}\mu(x)\bigg)^{1/p'}.
\]
As $p\to1^+$, we have $q\to q_1$ and $p'\to\infty$. Fix $R\in\mathcal R$. Since $\mu(R)<\infty$, the normalized measure $\mu_R=\mu|_R/\mu(R)$ is a probability measure on $R$. By the standard properties of $L^p$ norms on a finite measure space, we have
\[
\lim_{p'\to\infty}
\bigg(\frac{1}{\mu(R)}\sum_{x\in R}\omega(x)^{-p'}\mu(x)\bigg)^{1/p'}
=
\operatorname*{ess\,sup}_{x\in R}\omega(x)^{-1},
\]
and, since $q\to q_1$ as $p\to1^+$,
\[
\lim_{q\to q_1}
\bigg(\frac{1}{\mu(R)}\sum_{x\in R}\omega(x)^q\mu(x)\bigg)^{1/q}
=
\bigg(\frac{1}{\mu(R)}\sum_{x\in R}\omega(x)^{q_1}\mu(x)\bigg)^{1/q_1}.
\]
Hence $F_R(p)\to F_R(1)$, where
\[
F_R(1)
:=
\bigg(\frac{1}{\mu(R)}\sum_{x\in R}\omega(x)^{q_1}\mu(x)\bigg)^{1/q_1}
\operatorname*{ess\,sup}_{x\in R}\omega(x)^{-1}.
\]
Therefore, we have
\begin{align*}
	\liminf_{p\to1^+}[\omega]_{A(p,q)(\mu)}
	&=
	\liminf_{p\to1^+}\sup_{R\in\mathcal R}F_R(p)
	\ge
	\sup_{R\in\mathcal R}\lim_{p\to1^+}F_R(p)\\
	&=
	\sup_{R\in\mathcal R}F_R(1)
	=
	[\omega]_{A(1,q_1)(\mu)}.
\end{align*}

Combining the upper and lower bounds yields
\[
\lim_{p\to1^+}[\omega]_{A(p,q)(\mu)}
=
[\omega]_{A(1,q_1)(\mu)}.
\]

{\rm(v)} For $p>1$, property (ii) states that if $\omega\in A(p,q)(\mu)$, then $\omega^q\in A_q(\mu)$, and
\[
[\omega^q]_{A_q(\mu)}\lesssim [\omega]_{A(p,q)(\mu)}^q.
\]
Hence
\[
[\omega]_{A(p,q)(\mu)}
\gtrsim
[\omega^q]_{A_q(\mu)}^{1/q}.
\]

For $p=1$, property (iii) states that $\omega\in A(1,q)(\mu)$ if and only if $\omega^q\in A_1(\mu)$, and
\[
[\omega]_{A(1,q)(\mu)}
=
[\omega^q]_{A_1(\mu)}^{1/q}.
\]
The alternative characterizations of $A(1,q)(\mu)$ and $A(p,q)(\mu)$ in terms of suprema over $f$ now follow from the classical characterization of $A_r$ weights; see \cite[Proposition 7.1.5 (8)]{G2014}. The argument applies with cubes replaced by admissible trapezoids and the Lebesgue measure replaced by the flow measure $\mu$.

{\rm(vi)} For $1<p<1/\gamma$, property (i) gives
\[
\omega\in A(p,q)(\mu)
\quad\text{if and only if}\quad
\omega^q\in A_{1+q/p'}(\mu),
\]
and
\[
[\omega]_{A(p,q)(\mu)}
=
[\omega^q]_{A_{1+q/p'}(\mu)}^{1/q}.
\]
By the classical property of $A_r$ weights, see \cite[Proposition 7.1.5 (5)]{G2014},
\[
[\omega^q]_{A_{1+q/p'}(\mu)}\ge 1,
\]
with equality if and only if $\omega^q$ is constant. Hence
\[
[\omega]_{A(p,q)(\mu)}\ge 1,
\]
and equality holds if and only if $\omega$ is constant.

For $p=1$, property (iii) gives
\[
\omega\in A(1,q)(\mu)
\quad\text{if and only if}\quad
\omega^q\in A_1(\mu),
\]
and
\[
[\omega]_{A(1,q)(\mu)}
=
[\omega^q]_{A_1(\mu)}^{1/q}.
\]
Applying the same classical property with $r=1$ yields
\[
[\omega^q]_{A_1(\mu)}\ge 1,
\]
with equality if and only if $\omega^q$ is constant. Therefore
\[
[\omega]_{A(1,q)(\mu)}\ge 1,
\]
and equality holds if and only if $\omega$ is constant.

{\rm(vii)} For $p=1$, the definition of $A(1,q)(\mu)$ gives, for every $R\in\mathcal R$,
\[
\bigg(\frac{1}{\mu(R)}\sum_{x\in R}\omega(x)^q\mu(x)\bigg)^{1/q}
\operatorname*{ess\,sup}_{x\in R}\omega(x)^{-1}
\le [\omega]_{A(1,q)(\mu)}.
\]
Hence
\[
\bigg(\frac{1}{\mu(R)}\sum_{x\in R}\omega(x)^q\mu(x)\bigg)^{1/q}
\le [\omega]_{A(1,q)(\mu)}
\operatorname*{ess\,inf}_{x\in R}\omega(x).
\]
Taking the $q$-th power and multiplying by $\mu(R)$ yields
\[
\omega_{q,\mu}(R)
\le
[\omega]^q_{A(1,q)(\mu)}
\mu(R)
\Big(\operatorname*{ess\,inf}_{x\in R}\omega(x)\Big)^q.
\tag{1}
\]
Taking $R=R_1$ in (1) and using $\mu(R_1)\le C_0\mu(R_1\cap R_2)$, one obtains
\[
\omega_{q,\mu}(R_1)
\le
C_0 [\omega]^q_{A(1,q)(\mu)}
\mu(R_1\cap R_2)
\Big(\operatorname*{ess\,inf}_{x\in R_1}\omega(x)\Big)^q.
\]
Since $R_1\cap R_2\subset R_1$, it follows that
\[
\operatorname*{ess\,inf}_{x\in R_1}\omega(x)
\le
\operatorname*{ess\,inf}_{x\in R_1\cap R_2}\omega(x).
\]
Therefore
\begin{align*}
\omega_{q,\mu}(R_1)
&\le
C_0[\omega]^q_{A(1,q)(\mu)}
\mu(R_1\cap R_2)
\Big(\operatorname*{ess\,inf}_{x\in R_1\cap R_2}\omega(x)\Big)^q \\
&\le
C_0[\omega]^q_{A(1,q)(\mu)}
\sum_{x\in R_1\cap R_2}\omega(x)^q\mu(x) \\
&=
C_0[\omega]^q_{A(1,q)(\mu)}
\omega_{q,\mu}(R_1\cap R_2).
\end{align*}
Finally, since $R_1\cap R_2\subset R_2$, one has
\[
\omega_{q,\mu}(R_1\cap R_2)\le \omega_{q,\mu}(R_2).
\]
Thus
\[
\omega_{q,\mu}(R_1)
\leq C_0[\omega]^q_{A(1,q)(\mu)}
\omega_{q,\mu}(R_2).
\]

For $1<p<1/\gamma$, property (ii) gives $\omega^q\in A_q(\mu)$ and
\[
[\omega^q]_{A_q(\mu)}
\le
C_1(p,q)[\omega]^q_{A(p,q)(\mu)},
\]
where $C_1(p,q)$ depends on $p$ and $q$. For every $R\in\mathcal R$, the standard estimate for $A_q(\mu)$ weights yields
\[
\frac{1}{\mu(R)}\sum_{x\in R}\omega(x)^q\mu(x)
\le
[\omega^q]_{A_q(\mu)}
\operatorname*{ess\,inf}_{x\in R}\omega(x)^q.
\]
Multiplying by $\mu(R)$, one obtains
\[
\omega_{q,\mu}(R)
\le
C_1(p,q)[\omega]^q_{A(p,q)(\mu)}
\mu(R)
\Big(\operatorname*{ess\,inf}_{x\in R}\omega(x)\Big)^q.
\tag{2}
\]
Taking $R=R_1$ in (2) and using $\mu(R_1)\le C_0\mu(R_1\cap R_2)$, we get
\begin{align*}
\omega_{q,\mu}(R_1)
&\le
C_0 C_1(p,q)[\omega]^q_{A(p,q)(\mu)}
\mu(R_1\cap R_2)
\Big(\operatorname*{ess\,inf}_{x\in R_1\cap R_2}\omega(x)\Big)^q \\
&\le
C_0 C_1(p,q)[\omega]^q_{A(p,q)(\mu)}
\omega_{q,\mu}(R_1\cap R_2) \\
&\le
C_0 C_1(p,q)[\omega]^q_{A(p,q)(\mu)}
\omega_{q,\mu}(R_2).
\end{align*}
Setting $C_{p,q,C_0}=C_0 C_1(p,q)$ completes the proof.
\end{proof}

We now relate Definition \ref{dfn2-1} to the classical setting. For weights on $\mathbb Z$, the $A_p$ classes were introduced in \cite{C1976}; in the tree setting, the analogue of these classes was studied in \cite{AFM2024}. There it is proved that a level weight $w=W\circ\ell$, that is,
\[
w(x)=W(\ell(x)),\qquad x\in V,
\]
for some weight $W$ on $\mathbb Z$, is an $A_p(\mu)$ weight on $T$ if and only if $W$ is an $A_p(\mathbb Z)$ weight, with equal constants:
\[
[w]_{A_p(\mu)}=[W]_{A_p(\mathbb Z)}.
\] Our aim is to extend this correspondence to the fractional range. To this end, we first recall the classical $A(p,q)(\mathbb Z)$ classes.

\begin{dfn}
For $1<p<\infty$ and $1\le q<\infty$, a weight $W:\mathbb Z\to(0,\infty)$ is called an $A(p,q)(\mathbb Z)$ weight if
\[
[W]_{A(p,q)(\mathbb Z)}
:=
\sup_{I\ \text{interval in}\ \mathbb Z}
\bigg(\frac{1}{|I|}\sum_{x\in I}W(x)^q\bigg)^{1/q}
\bigg(\frac{1}{|I|}\sum_{x\in I}W(x)^{-p'}\bigg)^{1/p'}
<\infty,
\]
where $|\cdot|$ denotes the counting measure on $\mathbb Z$. Similarly, $W$ is called an $A(1,q)(\mathbb Z)$ weight if
\[
[W]_{A(1,q)(\mathbb Z)}
:=
\sup_{I\ \text{interval in}\ \mathbb Z}
\bigg(\frac{1}{|I|}\sum_{x\in I}W(x)^q\bigg)^{1/q}
\|W^{-1}\|_{\ell^\infty(I)}
<\infty.
\]
\end{dfn}

Let $\mathcal L(T)$ denote the class of level weights on $T$, namely the weights that are constant along level sets. Equivalently, $\omega\in\mathcal L(T)$ if there exists a positive function $W:\mathbb Z\to(0,\infty)$ such that $\omega=W\circ\ell$. Such weights are natural in this setting, since the flow measure $\mu$ assigns the same mass to each level inside an admissible trapezoid, so that the $A(p,q)(\mu)$ condition for level weights reduces to a condition on intervals in $\mathbb Z$.

The following result shows that, for such weights, Definition \ref{dfn2-1} reproduces the classical $A(p,q)$ condition exactly. This extends the $A_p$ case treated in \cite[Theorem 2.4]{AFM2024} to the fractional range.
\begin{prop}
Let $0<\gamma<1$ and $1\le p<1/\gamma$, and let $q\in[1,\infty)$ satisfy ${1}/{p}-{1}/{q}=\gamma$. For every level weight $\omega=W\circ\ell\in\mathcal L(T)$, the following hold:
\begin{enumerate}
	\item[(i)] If $1<p<1/\gamma$, then $\omega\in A(p,q)(\mu)$ if and only if $W\in A(p,q)(\mathbb Z)$, and
	\[
	[\omega]_{A(p,q)(\mu)}=[W]_{A(p,q)(\mathbb Z)}.
	\]
	\item[(ii)] If $p=1$, then $\omega\in A(1,q)(\mu)$ if and only if $W\in A(1,q)(\mathbb Z)$, and the corresponding constants are equal.
\end{enumerate}
\end{prop}
The case $p=q$ corresponds to the limiting case $\gamma=0$ and is covered by \cite[Theorem 2.4]{AFM2024}.

\begin{proof}
Assume first that $1<p<1/\gamma$. For every $R=R_{h_1}^{h_2}(x_0)\in\mathcal R$, since $\mu$ is a flow measure, we have
$$
\mu(R_{h_1}^{h_2}(x_0))=(h_2-h_1)\mu(x_0),\qquad x_0\in V.
$$ See Figure \ref{fig:trapezoid-level} for the level decomposition.
\begin{figure}[htbp]
	\centering
\def\sc{1.00}    
\pgfmathsetlengthmacro\lwTree {0.38pt*\sc}   
\pgfmathsetlengthmacro\lwTrap {0.95pt*\sc}   
\pgfmathsetlengthmacro\aOn    {3.2pt*\sc}    
\pgfmathsetlengthmacro\aOff   {2.4pt*\sc}
\pgfmathsetlengthmacro\lOn    {2.4pt*\sc}    
\pgfmathsetlengthmacro\lOff   {2.0pt*\sc}
\pgfmathsetlengthmacro\arA    {3.0pt*\sc}    
\pgfmathsetlengthmacro\arB    {2.0pt*\sc}
\pgfmathsetmacro\fsz{9*\sc}
\pgfmathsetmacro\fsl{11*\sc}
\begin{tikzpicture}[scale=\sc, transform shape,
	line cap=round, line join=round,
	tree/.style = {draw=black!80, line width=\lwTree},
	lev/.style  = {draw=black!70, line width=\lwTree, dash pattern=on \lOn off \lOff},
	trap/.style = {draw=black, line width=\lwTrap, fill=none},   
	vert/.style = {circle, fill=black!80, draw=none, inner sep=0pt, minimum size=3.0pt},
	redd/.style = {circle, fill=red,      draw=none, inner sep=0pt, minimum size=3.4pt},
	axis/.style = {draw=black!65, line width=\lwTree, dash pattern=on \aOn off \aOff},
	lab/.style  = {inner sep=1.5pt, font=\fontsize{\fsz}{\fsl}\selectfont},
	]
	\def\hs{0.55}    
	\def\dy{0.80}    
	\def\yXi{1.00}   
	\def\yTA{-1.36}\def\wTA{1.85}  
	\def\yTB{-4.44}\def\wTB{2.75}  
	\def\yTop{-1.60} 
	\def\yMid{-2.40} 
	\def\yPre{-3.40} 
	\def\yBot{-4.20} 
	\def\extT{2.15}  
	\def\extB{3.00}  
	\def\xpA{0.70}\def\xpB{2.10}   
	\def\offc{0.35}                
	
	\draw[axis] (0,\yXi) -- (0,0);
	\node[lab, anchor=south] at (0,\yXi+0.05) {$\xi$};
	\node[lab] at (0,\yXi+0.62) {root at infinity};
	
	\foreach \k/\nb in {0/1,1/2,2/4}{
		\pgfmathtruncatemacro{\nbm}{\nb-1}
		\foreach \j in {0,...,\nbm}{
			\draw[tree]
			({(2*\j+1)*2^(2-\k)*\hs-4*\hs},{-\k*\dy})
			--({(4*\j+1)*2^(1-\k)*\hs-4*\hs},{-(\k+1)*\dy});
			\draw[tree]
			({(2*\j+1)*2^(2-\k)*\hs-4*\hs},{-\k*\dy})
			--({(4*\j+3)*2^(1-\k)*\hs-4*\hs},{-(\k+1)*\dy});
		}
	}
	\foreach \k/\nb in {1/2,2/4,3/8}{
		\pgfmathtruncatemacro{\nbm}{\nb-1}
		\foreach \j in {0,...,\nbm}{
			\ifnum\k=1
			\node[vert] at ({(2*\j+1)*2^(2-\k)*\hs-4*\hs},{-\k*\dy}) {};
			\else
			\node[redd] at ({(2*\j+1)*2^(2-\k)*\hs-4*\hs},{-\k*\dy}) {};
			\fi
		}
	}
	
	\draw[trap] (-\wTA,\yTA) -- (\wTA,\yTA) -- (\wTB,\yTB) -- (-\wTB,\yTB) -- cycle;
	
	\draw[lev] (-\extT,\yTop) -- (\extT,\yTop);
	\node[lab, anchor=east] at (-\extT-0.08,\yTop) {$W(h)$};
	\node[lab, anchor=west] at ( \extT+0.10,\yTop) {$\mu(x_0)$};
	\node[lab, anchor=west] at ( \extT+1.18,\yTop) {$l(x_0)-h_1$};
	\draw[lev] (-\extT,\yMid) -- (\extT,\yMid);
	\node[lab, anchor=west] at ( \extT+0.12,\yMid) {$l(x)=h$};
	
	\foreach \xx in {-1.40,0,1.40}{ \node[lab] at (\xx,-2.92) {$\vdots$}; }
	
	\draw[lev] (-\xpB-\offc,\yPre) -- (\xpB+\offc,\yPre);
	\foreach \xx in {-\xpB,-\xpA,\xpA,\xpB}{ \node[redd] at (\xx,\yPre) {}; }
	\draw[lev] (-\extB,\yBot) -- (\extB,\yBot);
	\pgfmathsetmacro\cmA{\xpA-\offc}
	\pgfmathsetmacro\cpA{\xpA+\offc}
	\pgfmathsetmacro\cmB{\xpB-\offc}
	\pgfmathsetmacro\cpB{\xpB+\offc}
	\foreach \s in {1,-1}{
		\draw[tree] (\s*\xpA,\yPre) -- (\s*\cmA,\yBot);
		\draw[tree] (\s*\xpA,\yPre) -- (\s*\cpA,\yBot);
		\draw[tree] (\s*\xpB,\yPre) -- (\s*\cmB,\yBot);
		\draw[tree] (\s*\xpB,\yPre) -- (\s*\cpB,\yBot);
	}
	\foreach \xx in {-\cmB,-\cmA,\cpA,\cpB,-\cpA,-\cpB,\cmA,\cmB}{ \node[redd] at (\xx,\yBot) {}; }
	\node[lab, anchor=east] at (-\extB-0.08,\yBot) {$W(h)$};
	\node[lab, anchor=west] at ( \extB+0.10,\yBot) {$\mu(x_0)$};
	\node[lab, anchor=west] at ( \extB+1.18,\yBot) {$l(x_0)-h_2+1$};
	
	\node[lab] at (0,\yTB-0.62) {$R^{h_2}_{h_1}(x_0)$};
	\node[vert] at (0,0) {};
	\node[lab, anchor=east] at (-0.20,0) {$x_0$};
	
	
\end{tikzpicture}
\caption{Decomposition of the trapezoid \(R_{h_1}^{h_2}(x_0)\) by levels. Each level \(l(x)=h\) has total flow measure \(\mu(x_0)\), and \(\omega(x)=W(h)\) is constant on that level.}
\label{fig:trapezoid-level}
\end{figure}
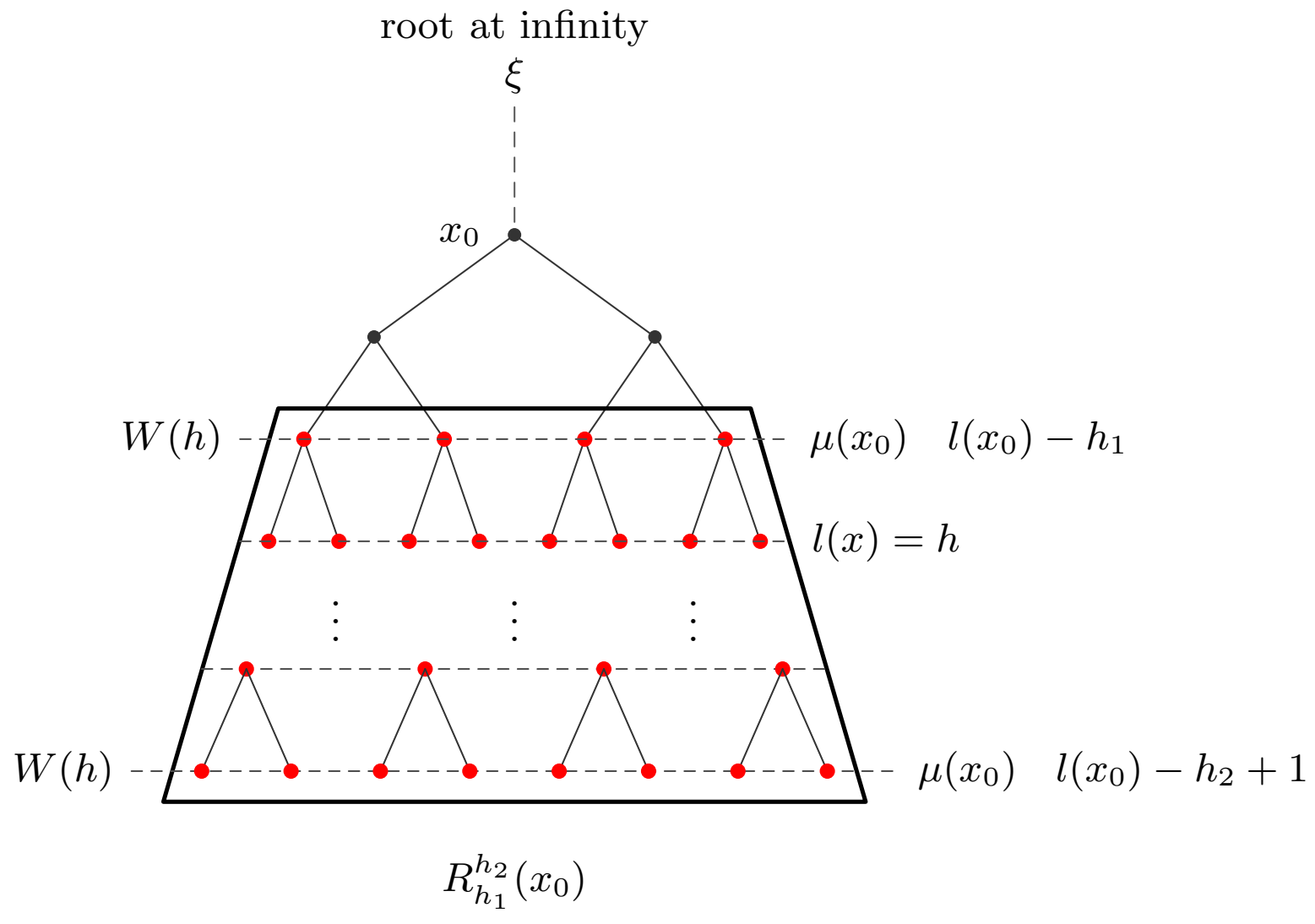

Thus we have
\begin{align*}
	F_R(p)
	&:=\Big(\frac{1}{\mu(R)}\sum_{x\in R}\omega(x)^q\mu(x)\Big)^{1/q}
	\bigg(\frac{1}{\mu(R)}\sum_{x\in R}\omega(x)^{-p'}\mu(x)\bigg)^{1/p'} \\
	&=\Bigg(\frac{1}{\mu(x_0)(h_2-h_1)}
	\sum_{h=l(x_0)-h_2+1}^{l(x_0)-h_1}
	W(h)^q\sum_{\substack{x\le x_0\\ l(x)=h}}\mu(x)\Bigg)^{1/q} \\
	&\quad\times
	\Bigg(\frac{1}{\mu(x_0)(h_2-h_1)}
	\sum_{h=l(x_0)-h_2+1}^{l(x_0)-h_1}
	W(h)^{-p'}\sum_{\substack{x\le x_0\\ l(x)=h}}\mu(x)\Bigg)^{1/p'}.
\end{align*}
Since $\mu$ is a flow measure,
$$
\sum_{\substack{x\le x_0\\ l(x)=h}}\mu(x)=\mu(x_0)
\qquad\text{for every }h\le l(x_0).
$$
Hence
$$
F_R(p)=
\Bigg(\frac{1}{h_2-h_1}\sum_{h=l(x_0)-h_2+1}^{l(x_0)-h_1}W(h)^q\Bigg)^{1/q}
\Bigg(\frac{1}{h_2-h_1}\sum_{h=l(x_0)-h_2+1}^{l(x_0)-h_1}W(h)^{-p'}\Bigg)^{1/p'}.
$$

The interval
\[
[l(x_0)-h_2+1,\ l(x_0)-h_1]\subset\mathbb Z
\]
has length $h_2-h_1$, and the admissibility condition $2\le h_2/h_1\le\beta$ is independent of the location of the interval. Therefore,
\[
[\omega]_{A(p,q)(\mu)}
\le [W]_{A(p,q)(\mathbb Z)}.
\]

Conversely, let $I=\{a,a+1,\dots,b\}\subset\mathbb Z$ be an arbitrary interval. Set
\[
h_1=b-a+1,\qquad h_2=2h_1,\qquad l(x_0)=2b-a+1.
\]
Then $2\le h_2/h_1=2\le\beta$, and
\[
l(x_0)-h_1=b,\qquad l(x_0)-h_2+1=a.
\]
Since the level function on an infinite tree takes every integer value, there exists $x_0\in V$ with $l(x_0)=2b-a+1$. For such $x_0$, the trapezoid $R_{h_1}^{h_2}(x_0)$ is admissible, and the computation above gives
\[
F_{R_{h_1}^{h_2}(x_0)}(p)
=
\bigg(\frac{1}{|I|}\sum_{h\in I}W(h)^q\bigg)^{1/q}
\bigg(\frac{1}{|I|}\sum_{h\in I}W(h)^{-p'}\bigg)^{1/p'}.
\]
Taking the supremum over all intervals $I\subset\mathbb Z$ yields
\[
[\omega]_{A(p,q)(\mu)}
\ge [W]_{A(p,q)(\mathbb Z)}.
\]
Thus equality holds for $1<p<1/\gamma$.

Finally, consider the case $p=1$. Then $q=1/(1-\gamma)$. By property (iv) of Proposition \ref{prop 2-1}, if $\omega\in A(1,q)(\mu)$, then
\[
\lim_{p\to1^+}[\omega]_{A(p,q)(\mu)}=[\omega]_{A(1,q)(\mu)}.
\]
On the other hand, for every $p>1$ we have just proved
\[
[\omega]_{A(p,q)(\mu)}=[W]_{A(p,q)(\mathbb Z)}.
\]
Passing to the limit as $p\to1^+$ and using the analogous limit property for $A(p,q)(\mathbb Z)$ weights, we obtain
\[
[\omega]_{A(1,q)(\mu)}=[W]_{A(1,q)(\mathbb Z)}.
\]
The equivalence $\omega\in A(1,q)(\mu)$ if and only if $W\in A(1,q)(\mathbb Z)$ follows from the same limiting argument. This completes the proof.
\end{proof}

The proof relies on the following covering lemma from \cite{AFM2024}, together with the new weight class $A(p,q)(\mu)$ and the weighted maximal estimates developed below.
\begin{lem}\cite{AFM2024}
Let $T$ be an infinite tree rooted at infinity, $\mu$ a flow measure on $T$, and let $\beta\ge 12$ be fixed. Denote by $\mathcal{R}$ the family of all admissible trapezoids and singletons. Let
\[
R=R_{h_1}^{h_2}(x)\in\mathcal{R},
\qquad h_1,h_2\in\mathbb{N},\quad h_1<h_2,\quad 2\le \frac{h_2}{h_1}\le \beta,
\]
and let
\[
\widetilde{R}:=R_{\lceil h_1/\beta\rceil}^{h_2\beta}(x)
\]
be its envelope. Define
\begin{align*}
R_0 &:= R_{\lceil h_1/\beta\rceil}^{h_1}(x), & \bar{R}_0 &:= R_{1}^{\beta}(x),\\
R_1 &:= R_{\lceil (h_1+h_2)/(2\beta)\rceil}^{\lfloor (h_1+h_2)/2\rfloor}(x), & \bar{R}_1 &:= R_{\lceil \beta/2\rceil}^{\lfloor \beta/2\rfloor\beta}(x),\\
R_2 &:= R_{\lfloor (h_1+h_2)/2\rfloor}^{\lfloor (h_1+h_2)/2\rfloor\beta}(x), & R_3 &:= R_{\lceil (h_1+h_2)/2\rceil\lfloor \beta/2\rfloor}^{\lceil (h_1+h_2)/2\rceil\lfloor \beta/2\rfloor\beta}(x).
\end{align*}
Here and below, $A\approx_{\beta} B$ means that there exists a constant $C\ge 1$, depending only on $\beta$, such that
\[
C^{-1}A\le B\le CA.
\]
Then the following assertions hold.

{\rm(i)} If $h_1\ge 3$, then
\[
\widetilde{R}\subset R_0\cup R_1\cup R_2\cup R_3.
\]
Moreover, $R_j\in\mathcal{R}$ and $\mu(R)\approx_{\beta}\mu(R_j)$ for $j=0,1,2,3$, and
\[
\mu(R)\approx_{\beta}\mu(R_0\cap R_1)
\approx_{\beta}\mu(R_1\cap R)
\approx_{\beta}\mu(R_2\cap R)
\approx_{\beta}\mu(R_3\cap R_2).
\]

{\rm(ii)} If $h_1=1$ and $h_2\ge 3$, or if $h_1=2$, then
\[
\widetilde{R}\subset \bar{R}_0\cup \bar{R}_1\cup R_2\cup R_3,
\]
and
\[
\mu(R)\approx_{\beta}\mu(\bar{R}_0\cap R_1)
\approx_{\beta}\mu(\bar{R}_1\cap R_2)
\approx_{\beta}\mu(R_2\cap R)
\approx_{\beta}\mu(R_3\cap R_2)
\approx_{\beta}\mu(R_3).
\]

{\rm(iii)} If $h_1=1$ and $h_2=2$, then
\[
\widetilde{R}\subset\bar{R}_0\cup\bar{R}_1,
\]
and
\[
\mu(R)=\mu(\bar{R}_0\cap R)
\approx_{\beta}\mu(\bar{R}_0\cap \bar{R}_1)
\approx_{\beta}\mu(\bar{R}_0)
\approx_{\beta}\mu(\bar{R}_1).
\]
\end{lem}

\begin{prop}\label{prop 2-3}
Let $0<\gamma<1$ and $1\le p<1/\gamma$, and let $q$ be defined by 
${1}/{p}-{1}/{q}=\gamma$. Then there exists a constant $C_{\beta,p,q}>0$ such that, for every $\omega\in A(p,q)(\mu)$ (with $A(p,q)(\mu)$ replaced by $A(1,q)(\mu)$ when $p=1$),
\[
\omega_{q,\mu}(\widetilde{R})\le C_{\beta,p,q}[\omega]_{A(p,q)(\mu)}^{2q}\omega_{q,\mu}(R)
\]
for every $R\in\mathcal R$.
\end{prop}

\begin{proof}
We split the proof into two cases.

\medskip
\textbf{Case 1: $p>1$.}
By Proposition \ref{prop 2-1}(ii), the assumption $\omega\in A(p,q)(\mu)$ implies $\omega^q\in A_q(\mu)$, and there exists a constant $C_1>0$, depending only on $p$ and $q$, such that
\[
[\omega^q]_{A_q(\mu)}\le C_1(p,q) [\omega]_{A(p,q)(\mu)}^q.
\]
Apply \cite[Theorem 2.7]{AFM2024} to the weight $\widetilde{\omega}=\omega^q$ and the exponent $q$. Since
\[
\widetilde{\omega}_\mu(R)=\sum_{x\in R}\omega(x)^q\mu(x)=\omega_{q,\mu}(R),
\]
we obtain a constant $C_{\beta,q}>0$ such that
\[
\omega_{q,\mu}(\widetilde R)
= \widetilde{\omega}_\mu(\widetilde R)
\le C_{\beta,q} [\widetilde{\omega}]_{A_q(\mu)}^2 \widetilde{\omega}_\mu(R)
= C_{\beta,q} [\omega^q]_{A_q(\mu)}^2 \omega_{q,\mu}(R).
\]
Combining the two estimates gives
\[
\omega_{q,\mu}(\widetilde R)
\le C_{\beta,q} (C_1(p,q))^2 [\omega]_{A(p,q)(\mu)}^{2q} \omega_{q,\mu}(R).
\]
Setting $C_{\beta,p,q}=C_{\beta,q} (C_1(p,q))^2$ yields the desired inequality.

\medskip
\textbf{Case 2: $p=1$.}
By Proposition \ref{prop 2-1}(iii), $\omega\in A(1,q)(\mu)$ if and only if $\omega^q\in A_1(\mu)$, and
\[
[\omega]_{A(1,q)(\mu)}=[\omega^q]_{A_1(\mu)}^{1/q}.
\]
Applying \cite[Theorem 2.7]{AFM2024} with $\widetilde{\omega}=\omega^q$ and $p=1$ (i.e., $A_1(\mu)$), we get
\[
\omega_{q,\mu}(\widetilde R)
\le C_{\beta,1} [\omega^q]_{A_1(\mu)}^2 \omega_{q,\mu}(R)
= C_{\beta,1} [\omega]_{A(1,q)(\mu)}^{2q} \omega_{q,\mu}(R).
\]
Thus the desired estimate holds. This completes the proof.
\end{proof}

\section{Characterization of Weighted Boundedness for the Fractional Maximal Operator on Trees}
The main result of this section is the following characterization of the weighted boundedness of the fractional maximal operator $M_{\mu,\gamma}$ on trees.

\begin{thm}\label{thm 3-1}
Let $0\le \gamma<1$ and $1\le p<1/\gamma$, and let $q\in[1,\infty)$ satisfy ${1}/{p}-{1}/{q}=\gamma$. Then the following hold.

{\rm(i)} If $1<p<1/\gamma$, then $M_{\mu,\gamma}$ is bounded from $\ell^p_{\mu}(\omega^p)$ to $\ell^q_{\mu}(\omega^q)$ if and only if $\omega\in A(p,q)(\mu)$.

{\rm(ii)} If $p=1$, then $M_{\mu,\gamma}$ is bounded from $\ell^1_{\mu}(\omega)$ to $\ell^{q,\infty}_{\mu}(\omega^q)$ if and only if $\omega\in A(1,q)(\mu)$.
\end{thm}

In the case $\gamma=0$, the result reduces to the corresponding conclusion established in \cite{AFM2024}.

The proof of Theorem \ref{thm 3-1} is divided into two parts: necessity and sufficiency. We first establish the necessity of the weight condition.

\begin{prop}\label{prop 3-1}
Let $0\le \gamma<1$ and $1\le p<1/\gamma$, and let $q\in[1,\infty)$ satisfy ${1}/{p}-{1}/{q}=\gamma$. Then the following hold.

{\rm(i)} If $1<p<1/\gamma$ and $M_{\mu,\gamma}$ is bounded from $\ell^p_{\mu}(\omega^p)$ to $\ell^q_{\mu}(\omega^q)$, then $\omega\in A(p,q)(\mu)$.

{\rm(ii)} If $p=1$ and $M_{\mu,\gamma}$ is bounded from $\ell^1_{\mu}(\omega)$ to $\ell^{q,\infty}_{\mu}(\omega^q)$, then $\omega\in A(1,q)(\mu)$.
\end{prop}

\begin{proof}[Proof of Proposition \ref{prop 3-1}]
We first prove (i). Assume that $1<p<1/\gamma$ and that $M_{\mu,\gamma}$ is bounded from $\ell^p_{\mu}(\omega^p)$ to $\ell^q_{\mu}(\omega^q)$, where ${1}/{p}-{1}/{q}=\gamma$.
Thus there exists a constant $C>0$ such that
\[
\bigg(\sum_{x\in V}\big(M_{\mu,\gamma}f(x)\omega(x)\big)^q\mu(x)\bigg)^{1/q}
\le C
\bigg(\sum_{x\in V}|f(x)|^p\omega(x)^p\mu(x)\bigg)^{1/p}
\]
for every $f$.

Fix $R\in\mathcal R$ and take
\[
f(x)=\omega(x)^{-p'}\chi_R(x).
\]
Then we have
\[
\bigg(\sum_{x\in V}|f(x)|^p\omega(x)^p\mu(x)\bigg)^{1/p}
=
\bigg(\sum_{x\in R}\omega(x)^{-p'}\mu(x)\bigg)^{1/p}.
\]

On the other hand, for every $x\in R$,
\[
M_{\mu,\gamma}f(x)
\ge
\frac{1}{\mu(R)^{1-\gamma}}\sum_{y\in R}f(y)\mu(y)
=
\frac{1}{\mu(R)^{1-\gamma}}\sum_{y\in R}\omega(y)^{-p'}\mu(y).
\]
Therefore
\[
\sum_{x\in R}\big(M_{\mu,\gamma}f(x)\omega(x)\big)^q\mu(x)
\ge
\bigg(\sum_{x\in R}\omega(x)^q\mu(x)\bigg)
\bigg(\frac{1}{\mu(R)^{1-\gamma}}\sum_{y\in R}\omega(y)^{-p'}\mu(y)\bigg)^q.
\]
Combining this with the boundedness assumption gives
\[
\bigg(\sum_{x\in R}\omega(x)^q\mu(x)\bigg)^{1/q}
\bigg(\frac{1}{\mu(R)^{1-\gamma}}\sum_{x\in R}\omega(x)^{-p'}\mu(x)\bigg)
\le
C
\bigg(\sum_{x\in R}\omega(x)^{-p'}\mu(x)\bigg)^{1/p}.
\]

Since
\[
1-\gamma
=
\frac{1}{q}+\frac{1}{p'},
\]
we obtain
\[
\bigg(\frac{1}{\mu(R)}\sum_{x\in R}\omega(x)^q\mu(x)\bigg)^{1/q}
\bigg(\frac{1}{\mu(R)}\sum_{x\in R}\omega(x)^{-p'}\mu(x)\bigg)^{1/p'}
\le
C.
\]
Taking the supremum over all $R\in\mathcal R$ gives
\[
[\omega]_{A(p,q)(\mu)}\le C<\infty,
\]
so $\omega\in A(p,q)(\mu)$.

We now prove (ii). Assume that $p=1$ and that $M_{\mu,\gamma}$ is bounded from $\ell^1_{\mu}(\omega)$ to $\ell^{q,\infty}_{\mu}(\omega^q)$, where $q=(1-\gamma)^{-1}$. 
Thus there exists a constant $C>0$ such that
\[
\lambda
\Bigg(\sum_{\{x\in V:\,M_{\mu,\gamma}f(x)>\lambda\}}\omega(x)^q\mu(x)\Bigg)^{1/q}
\le
C
\sum_{x\in V}|f(x)|\omega(x)\mu(x)
\]
for every $f\in \ell^1_{\mu}(\omega)$ and every $\lambda>0$.

Fix $R_1,R_2\in\mathcal R$ with $R_1\subset R_2$, and take $f=\chi_{R_1}$. For every $x\in R_2$, since $R_2$ is an admissible trapezoid containing $x$, we have
\[
M_{\mu,\gamma}f(x)
\ge
\frac{1}{\mu(R_2)^{1-\gamma}}\sum_{y\in R_2}f(y)\mu(y)
=
\frac{\mu(R_1)}{\mu(R_2)^{1-\gamma}}.
\]
Hence
\[
R_2\subset
\left\{
x\in V:\,
M_{\mu,\gamma}f(x)\ge
\frac{\mu(R_1)}{\mu(R_2)^{1-\gamma}}
\right\}.
\]
Set
\[
\lambda:=\frac{\mu(R_1)}{\mu(R_2)^{1-\gamma}}.
\]
Applying the weak-type estimate with this $\lambda$, we obtain
\[
\lambda
\bigg(\sum_{x\in R_2}\omega(x)^q\mu(x)\bigg)^{1/q}
\le
C
\sum_{x\in R_1}\omega(x)\mu(x).
\]
Substituting the value of $\lambda$ yields
\[
\frac{\mu(R_1)}{\mu(R_2)^{1-\gamma}}
\bigg(\sum_{x\in R_2}\omega(x)^q\mu(x)\bigg)^{1/q}
\le
C
\sum_{x\in R_1}\omega(x)\mu(x).
\]
Since $q=1/(1-\gamma)$, we have $(1-\gamma)q=1$ and hence
\[
\frac{1}{\mu(R_2)}\sum_{x\in R_2}\omega(x)^q\mu(x)
\le
C^q
\bigg(
\frac{1}{\mu(R_1)}\sum_{x\in R_1}\omega(x)\mu(x)
\bigg)^q.
\]

Now take $R_1=\{x_0\}$ for an arbitrary vertex $x_0\in V$. Then the above inequality becomes
\[
\bigg(
\frac{1}{\mu(R_2)}\sum_{x\in R_2}\omega(x)^q\mu(x)
\bigg)^{1/q}
\le
C\omega(x_0)
\]
for every admissible trapezoid $R_2$ containing $x_0$. In particular, for any fixed $R\in\mathcal R$ and any $x_0\in R$, taking $R_2=R$ gives
\[
\bigg(
\frac{1}{\mu(R)}\sum_{x\in R}\omega(x)^q\mu(x)
\bigg)^{1/q}
\le
C\omega(x_0).
\]
Since this holds for every $x_0\in R$, we obtain
\[
\bigg(
\frac{1}{\mu(R)}\sum_{x\in R}\omega(x)^q\mu(x)
\bigg)^{1/q}
\le
C\inf_{x_0\in R}\omega(x_0)
=
C\Big(\operatorname*{ess\,sup}_{x\in R}\omega(x)^{-1}\Big)^{-1}.
\]
Hence
\[
\bigg(
\frac{1}{\mu(R)}\sum_{x\in R}\omega(x)^q\mu(x)
\bigg)^{1/q}
\operatorname*{ess\,sup}_{x\in R}\omega(x)^{-1}
\le C.
\]
Taking the supremum over all $R\in\mathcal R$, we conclude that
\[
\sup_{R\in\mathcal R}
\bigg(
\frac{1}{\mu(R)}\sum_{x\in R}\omega(x)^q\mu(x)
\bigg)^{1/q}
\operatorname*{ess\,sup}_{x\in R}\omega(x)^{-1}
\le C,
\]
which is exactly the condition $\omega\in A(1,q)(\mu)$. This completes the proof of the necessity.
\end{proof}

For the covering argument of the level sets of the fractional maximal operator $M_{\mu,\gamma}$, the following geometric lemma from \cite{MFAM2023} is used. It compares two intersecting admissible trapezoids in terms of the mass of their roots.

\begin{lem}\label{lem 3-1}\cite{MFAM2023}
Let $R_1, R_2 \in \mathcal{R}$ with roots $x_1, x_2$, respectively. If
\[
R_1 \cap R_2 \neq \emptyset
\qquad\text{and}\qquad
\mu(x_1) \ge \mu(x_2),
\]
then
\[
R_2 \subset \widetilde{R}_1.
\]
\end{lem}

Since $\mu$ is not assumed to be doubling, the usual Vitali-type covering is not available. The following lemma provides a substitute based on an ordering of admissible trapezoids by root measures and on their envelopes.

\begin{lem}\label{lem 3-2}
Let $\lambda>0$ and $0\le\gamma<1$. For a function $f$ on $V$, set
\[
E_\lambda
=
\left\{
x \in V : M_{\mu,\gamma}f(x) > \lambda
\right\}.
\]
Then there exists a family $\{R_i\}_{i\in I}\subset \mathcal{R}$ such that
\[
E_\lambda = \bigcup_{i\in I} R_i.
\]
Moreover, there exists a pairwise disjoint subfamily $\{R_{i_k}\}_{k\in K}$ such that
\[
E_\lambda \subset \bigcup_{k\in K} \widetilde{R}_{i_k}.
\]
\end{lem}

\begin{proof}[Proof of Lemma \ref{lem 3-2}] Let $\lambda>0$, $0\le\gamma<1$, and assume $f\in\ell^1_\mu$. Let
\[
\mathcal F_0=
\Big\{
R\in\mathcal R:
\frac{1}{\mu(R)^{1-\gamma}}\sum_{y\in R}|f(y)|\mu(y)>\lambda
\Big\}.
\]
Then $E_\lambda=\bigcup_{R\in\mathcal F_0}R$. Indeed, the inclusion $E_\lambda\subset\bigcup_{R\in\mathcal F_0}R$ follows from the definition of $M_{\mu,\gamma}$, while the reverse inclusion is immediate.

For every $R\in\mathcal F_0$, we have
\[
\lambda\,\mu(R)^{1-\gamma}
<
\sum_{y\in R}|f(y)|\mu(y)
\le
\|f\|_{\ell^1_\mu}.
\]
Hence
\[
\mu(R)
<
\big({\|f\|_{\ell^1_\mu}}/{\lambda}\big)^{1/(1-\gamma)},
\]
so $\{\mu(R):R\in\mathcal F_0\}$ is bounded above.

Since $V$ is countable, we may enumerate $\mathcal F_0=\{R_1,R_2,\dots\}$ in such a way that
\[
\mu(x_{R_1})\ge \mu(x_{R_2})\ge\cdots.
\]
We now select a pairwise disjoint subfamily $\{R^*_k\}_{k\in K}$ by the following greedy procedure: set $R^*_1=R_1$, and after $R^*_1,\dots,R^*_k$ have been chosen, let $R^*_{k+1}$ be the first element of the enumeration that is disjoint from all previously chosen sets. If no such element exists, the procedure stops.

We claim that for every $R\in\mathcal F_0$ there exists $k\in K$ such that
\[
R\cap R^*_k\neq\emptyset
\qquad\text{and}\qquad
\mu(x_{R^*_k})\ge \mu(x_R).
\]

Assume, for contradiction, that $R$ is disjoint from every $R^*_k$. Write $R=R_m$. When the greedy procedure reaches $R_m$, it is disjoint from all previously chosen sets, so it would have been selected, a contradiction. Thus the claim holds. Moreover, by the ordering of the enumeration, $R^*_k$ appears no later than $R$, so $\mu(x_{R^*_k})\ge \mu(x_R)$.

Now fix $R\in\mathcal F_0$. By the claim, there exists $k\in K$ such that $R\cap R^*_k\neq\emptyset$ and $\mu(x_{R^*_k})\ge \mu(x_R)$. Applying Lemma \ref{lem 3-1} with $R_1=R^*_k$ and $R_2=R$, we obtain
\[
R\subset \widetilde{R^*_k}.
\]
Therefore
\[
E_\lambda=\bigcup_{R\in\mathcal F_0}R
\subset
\bigcup_{k\in K}\widetilde{R^*_k}.
\]
The subfamily $\{R^*_k\}_{k\in K}$ is pairwise disjoint by construction. This completes the proof.
\end{proof}

Having established the necessity of the $A(p,q)(\mu)$ condition in Proposition \ref{prop 3-1}, we now turn to the sufficiency. As is customary in weighted theory, the first step is to prove the weak-type endpoint estimate, which will serve as the foundation for the strong-type result. More precisely, we show that the $A(p,q)(\mu)$ condition already guarantees the boundedness of the fractional maximal operator $M_{\mu,\gamma}$ from $\ell^p_\mu(\omega^p)$ to the weak space $\ell^{q,\infty}_\mu(\omega^q)$. 

\begin{prop}
	Let $\omega$ be an $A(p,q)(\mu)$ weight with $0\le\gamma<1$, $1\le p<1/\gamma$, and ${1}/{p}-{1}/{q}=\gamma$.
	Then $M_{\mu,\gamma}$ is bounded from $\ell^p_\mu(\omega^p)$ to $\ell^{q,\infty}_\mu(\omega^q)$.
\end{prop}

\begin{proof}
We first consider the case $p>1$. Let $\lambda>0$ and set
\[
E_\lambda=\{x\in V:M_{\mu,\gamma}f(x)>\lambda\}.
\]
By Lemma \ref{lem 3-2}, there exists a pairwise disjoint subfamily $\{R_{j_k}\}_{k\in K}\subset\mathcal R$ such that
\[
E_\lambda\subset\bigcup_{k\in K}\widetilde{R}_{j_k}.
\]
Recall that the family $\mathcal F_0$ used in the proof of Lemma \ref{lem 3-2} is defined by
\[
\mathcal F_0=
\Big\{
R\in\mathcal R:
\frac{1}{\mu(R)^{1-\gamma}}\sum_{y\in R}|f(y)|\mu(y)>\lambda
\Big\}.
\]
For each $k\in K$, since $R_{j_k}\in\mathcal F_0$, we have
\[
\lambda < \frac{1}{\mu(R_{j_k})^{1-\gamma}}\sum_{y\in R_{j_k}}|f(y)|\mu(y).
\]
Hence
\begin{equation}\label{eq 3-1}
	\lambda^p \mu(R_{j_k})^{p(1-\gamma)}
	< \bigg(\sum_{y\in R_{j_k}}|f(y)|\mu(y)\bigg)^p.
\end{equation}

For each $k\in K$, Proposition \ref{prop 2-3} gives
\[
\omega_{q,\mu}(\widetilde{R}_{j_k})
\le C_{\beta,p,q}[\omega]_{A(p,q)(\mu)}^{2q}\omega_{q,\mu}(R_{j_k}).
\]
Raising both sides to the positive power $p/q$, we obtain
\[
\omega_{q,\mu}(\widetilde{R}_{j_k})^{p/q}
\le C' [\omega]_{A(p,q)(\mu)}^{2p}\omega_{q,\mu}(R_{j_k})^{p/q},
\]
where $C'=C_{\beta,p,q}^{p/q}$.

Since $\gamma\ge 0$ implies $p\le q$, we have $p/q\le 1$. Hence, by the subadditivity of $t\mapsto t^{p/q}$ on $[0,\infty)$, and using
\[
E_\lambda\subset\bigcup_{k\in K}\widetilde{R}_{j_k},
\]
we get
\begin{align*}
	\omega_{q,\mu}(E_\lambda)^{p/q}
	&\le \bigg(\sum_{k\in K}\omega_{q,\mu}(\widetilde{R}_{j_k})\bigg)^{p/q}
	\le \sum_{k\in K}\omega_{q,\mu}(\widetilde{R}_{j_k})^{p/q}\\
	&\le C' [\omega]_{A(p,q)(\mu)}^{2p}
	\sum_{k\in K}\omega_{q,\mu}(R_{j_k})^{p/q}.
\end{align*}

On the other hand, applying Hölder's inequality with exponents $p$ and $p'$,
\begin{equation}\label{eq 3-2}
	\bigg(\sum_{y\in R_{j_k}}|f(y)|\mu(y)\bigg)^p
	\le
	\bigg(\sum_{y\in R_{j_k}}|f(y)|^p\omega(y)^p\mu(y)\bigg)
	\bigg(\sum_{y\in R_{j_k}}\omega(y)^{-p'}\mu(y)\bigg)^{p/p'}.
\end{equation}
By the $A(p,q)(\mu)$ condition,
\[
\omega_{q,\mu}(R_{j_k})^{1/q}
\bigg(\sum_{y\in R_{j_k}}\omega(y)^{-p'}\mu(y)\bigg)^{1/p'}
\le
[\omega]_{A(p,q)(\mu)}\mu(R_{j_k})^{1-\gamma}.
\]
Raising this to the power $p$ and rearranging, we obtain
\begin{equation}\label{eq 3-3}
	\mu(R_{j_k})^{p(1-\gamma)}
	\ge
	\frac{\omega_{q,\mu}(R_{j_k})^{p/q}
		\left(\sum_{y\in R_{j_k}}\omega(y)^{-p'}\mu(y)\right)^{p/p'}}
	{[\omega]_{A(p,q)(\mu)}^p}.
\end{equation}

Note that $\omega>0$ and $R_{j_k}\neq\emptyset$, then
\[
\bigg(\sum_{y\in R_{j_k}}\omega(y)^{-p'}\mu(y)\bigg)^{p/p'}>0.
\]
Combining \eqref{eq 3-1}, \eqref{eq 3-2} and \eqref{eq 3-3}, we have 
\[
\lambda^p
\frac{\omega_{q,\mu}(R_{j_k})^{p/q}}{[\omega]_{A(p,q)(\mu)}^p}
\le
\sum_{y\in R_{j_k}}|f(y)|^p\omega(y)^p\mu(y).
\]
Thus
\[
\omega_{q,\mu}(R_{j_k})^{p/q}
\le
\frac{[\omega]_{A(p,q)(\mu)}^p}{\lambda^p}
\sum_{y\in R_{j_k}}|f(y)|^p\omega(y)^p\mu(y).
\]
Substituting this into the previous estimate for $\omega_{q,\mu}(E_\lambda)^{p/q}$, we obtain
\[
\omega_{q,\mu}(E_\lambda)^{p/q}
\le
C' [\omega]_{A(p,q)(\mu)}^{2p}
\frac{[\omega]_{A(p,q)(\mu)}^p}{\lambda^p}
\sum_{k\in K}\sum_{y\in R_{j_k}}|f(y)|^p\omega(y)^p\mu(y).
\]
Since the $R_{j_k}$ are pairwise disjoint,
\[
\sum_{k\in K}\sum_{y\in R_{j_k}}|f(y)|^p\omega(y)^p\mu(y)
\le
\|f\|_{\ell^p_\mu(\omega^p)}^p.
\]
Therefore
\[
\omega_{q,\mu}(E_\lambda)^{p/q}
\le
C'' \frac{[\omega]_{A(p,q)(\mu)}^{3p}}{\lambda^p}
\|f\|_{\ell^p_\mu(\omega^p)}^p,
\]
where $C''$ depends only on $\beta,p,q$. Raising both sides to the power $q/p$, we conclude that
\[
\omega_{q,\mu}(E_\lambda)
\le
C''' \frac{[\omega]_{A(p,q)(\mu)}^{3q}}{\lambda^q}
\|f\|_{\ell^p_\mu(\omega^p)}^q,
\]
which is the desired weak-type estimate.

The case $p=1$ is similar and even simpler. Indeed, when $p=1$ we have $q=1/(1-\gamma)$ and $\omega\in A(1,q)(\mu)$. By Proposition \ref{prop 2-1}(iii), $\omega^q\in A_1(\mu)$. Using the same covering argument and the analogue of Proposition \ref{prop 2-3} for $A_1(\mu)$ weights (or by a limiting argument from the $p>1$ case), one obtains the weak-type estimate from $\ell^1_\mu(\omega)$ to $\ell^{q,\infty}_\mu(\omega^q)$. We omit the details.
\end{proof}

Having established the weak-type estimate, we now turn to the strong-type boundedness. The following proposition shows that the $A(p,q)(\mu)$ condition is sufficient for $M_{\mu,\gamma}$ to be bounded from $\ell^p_\mu(\omega^p)$
to $\ell^q_\mu(\omega^q)$.

\begin{prop}
Let $0\le \gamma<1$ and $1<p<1/\gamma$, and let $q\in [p,\infty)$ be defined by ${1}/{p}-{1}/{q}=\gamma$. If $\omega\in A(p,q)(\mu)$, then $M_{\mu,\gamma}$ is bounded from
$\ell^p_\mu(\omega^p)$ to $\ell^q_\mu(\omega^q)$. More precisely, there exists
a constant $C=C(\beta,p,q)>0$ such that
\[
\|M_{\mu,\gamma}f\|_{\ell^q_\mu(\omega^q)}
\le
C[\omega]_{A(p,q)(\mu)}^{\alpha}
\|f\|_{\ell^p_\mu(\omega^p)}
\]
for every $f\in\ell^p_\mu(\omega^p)$, where $\alpha=\alpha(p,q)>0$ depends
only on $p$ and $q$.
\end{prop}

\begin{proof}
Let
\[
r=q(1-\gamma)=1+\frac{q}{p'}>1.
\]
Since $\omega\in A(p,q)(\mu)$, Proposition \ref{prop 2-1}(i) yields
$\omega^q\in A_r(\mu)$, and the definition gives
\[
[\omega^q]_{A_r(\mu)}=[\omega]_{A(p,q)(\mu)}^q.
\]
We shall use the following quantitative form of the weighted maximal
estimate from \cite[Theorem 2.9]{AFM2024}:
\[
\|M_\mu g\|_{\ell^r_\mu(\omega^q)}
\le C_{\beta,r}[\omega^q]_{A_r(\mu)}^{3/(r-1)}
\|g\|_{\ell^r_\mu(\omega^q)}.
\]
For completeness, we recall the dependence on the weight constant.
For a positive weight $a$, define the associated maximal operator by
\[
M^a_\mu g(x):=\sup_{R\in\mathcal R,\,x\in R}
\frac{\sum_{y\in R}|g(y)|a(y)\mu(y)}
{\sum_{y\in R}a(y)\mu(y)}.
\]
By \cite[Theorem 2.7 and Corollary 2.8]{AFM2024}, if
$a\in A_s(\mu)$, $s>1$, then the weak type $(1,1)$ constant of
$M^a_\mu$ with respect to $a\mu$ is bounded by
$C_{\beta,s}[a]_{A_s(\mu)}^2$. Interpolating with the trivial
$\ell^\infty$ bound yields, for every $t>1$,
\[
\|M^a_\mu\|_{\ell^t_\mu(a)\to\ell^t_\mu(a)}
\le C_{\beta,s,t}[a]_{A_s(\mu)}^{2/t}.
\]
Let $u=\omega^q$, $\sigma=u^{-1/(r-1)}$, and $K=[u]_{A_r(\mu)}$.
Then $\sigma\in A_{r'}(\mu)$ and
$[\sigma]_{A_{r'}(\mu)}=K^{1/(r-1)}$, and
\begin{align*}
	\|M_\mu\|_{\ell^r_\mu(u)\to\ell^r_\mu(u)}
	&\le K^{1/(r-1)}
	\|M^u_\mu\|_{\ell^{r'}_\mu(u)\to\ell^{r'}_\mu(u)}^{1/(r-1)}
	\|M^\sigma_\mu\|_{\ell^r_\mu(\sigma)\to\ell^r_\mu(\sigma)}\\
	&\le C_{\beta,r}K^{\frac{1}{r-1}+\frac{2}{r}+\frac{2}{r(r-1)}}
	=C_{\beta,r}K^{3/(r-1)}.
\end{align*}
This establishes the stated quantitative estimate.

If $\gamma=0$, then $q=r=p$ and $M_{\mu,\gamma}=M_\mu$, so the
desired estimate follows immediately with $\alpha=3p'$. We now assume
$0<\gamma<1$. Let $f\in\ell^p_\mu(\omega^p)$ and set
\[
g(x):=|f(x)|^{p/r}\omega(x)^{(p-q)/r},\qquad x\in V.
\]
Then
\[
\sum_{x\in V}|g(x)|^r\omega(x)^q\mu(x)
=\sum_{x\in V}|f(x)|^p\omega(x)^p\mu(x)<\infty,
\]
so $g\in\ell^r_\mu(\omega^q)$ and
\[
\|g\|_{\ell^r_\mu(\omega^q)}
=\|f\|_{\ell^p_\mu(\omega^p)}^{p/r}.
\]
Using the identities $p/q+p\gamma=1$ and $r=q(1-\gamma)$, we have
\[
|f(y)|=g(y)^{1-\gamma}
\bigl(|f(y)|^p\omega(y)^p\bigr)^\gamma.
\]
For every $R\in\mathcal R$, Hölder's inequality with exponents
$1/(1-\gamma)$ and $1/\gamma$ therefore gives
\begin{align*}
	\frac{1}{\mu(R)^{1-\gamma}}\sum_{y\in R}|f(y)|\mu(y)
	&\le
	\Big(\frac{1}{\mu(R)}\sum_{y\in R}g(y)\mu(y)\Big)^{1-\gamma}
	\Big(\sum_{y\in R}|f(y)|^p\omega(y)^p\mu(y)\Big)^\gamma\\
	&\le
	\Big(\frac{1}{\mu(R)}\sum_{y\in R}g(y)\mu(y)\Big)^{1-\gamma}
	\|f\|_{\ell^p_\mu(\omega^p)}^{p\gamma}.
\end{align*}
Taking the supremum over all $R\in\mathcal R$ containing $x$, we obtain
\[
M_{\mu,\gamma}f(x)
\le \bigl(M_\mu g(x)\bigr)^{1-\gamma}
\|f\|_{\ell^p_\mu(\omega^p)}^{p\gamma}.
\]
Since $q(1-\gamma)=r$, the weighted maximal estimate yields
\begin{align*}
	\|M_{\mu,\gamma}f\|_{\ell^q_\mu(\omega^q)}
	&\le
	\|M_\mu g\|_{\ell^r_\mu(\omega^q)}^{1-\gamma}
	\|f\|_{\ell^p_\mu(\omega^p)}^{p\gamma}\\
	&\le
	C_{\beta,p,q}[\omega^q]_{A_r(\mu)}^{\frac{3(1-\gamma)}{r-1}}
	\|g\|_{\ell^r_\mu(\omega^q)}^{1-\gamma}
	\|f\|_{\ell^p_\mu(\omega^p)}^{p\gamma}\\
	&=
	C_{\beta,p,q}[\omega]_{A(p,q)(\mu)}^{3p'(1-\gamma)}
	\|f\|_{\ell^p_\mu(\omega^p)},
\end{align*}
where in the last step we used $r-1=q/p'$ and
$p(1-\gamma)/r+p\gamma=p/q+p\gamma=1$. Thus the conclusion holds
with $\alpha=3p'(1-\gamma)>0$. This completes the proof.
\end{proof}

\section*{Declarations}

\subsection*{Conflict of Interest}
The authors declare no conflict of interest.

\subsection*{Funding}
The research was supported by National Natural Science Foundation of China (Grant No. 11971295).

\subsection*{AI Declaration}
The authors declare that no artificial intelligence tools were used in the preparation, generation, or editing of this manuscript. All mathematical arguments, proofs, and the writing of the text were carried out entirely by the authors.

\end{document}